\documentclass[12pt]{article}
\usepackage{amssymb,amscd,amsmath,amsthm}
\usepackage{latexsym,amstext}
\usepackage{latexsym,amstext}
\usepackage{color}
\usepackage{graphicx}
\usepackage{epstopdf}
\usepackage{cite}
\usepackage[utf8]{inputenc}

\begin{document}

\title{\bf A method to find approximate solutions of first order systems of non-linear ordinary equations.}

\author{J.J. Alvarez-S\'anchez$^{1}$, M. Gadella$^2$, L.P. Lara$^{3,4}$\\ \\
$^1$ Escuela de Inform\'atica de Segovia, Universidad de Valladolid, Segovia, Spain\\
$^2$ Departamento de F\'{\i}sica Te\'orica, At\'omica y Optica  and IMUVA, \\
Universidad de Va\-lladolid, 47011 Valladolid, Spain\\ 
$^3$ Instituto de F\'isica Rosario, CONICET-UNR, \\ 
Bv. 27 de Febrero, S2000EKF Rosario, Santa Fe,  Argentina.\\
$^4$ Departamento de Sistemas, Universidad del Centro Educativo\\
Latinoamericano, Av. Pellegrini 1332, S2000 Rosario, Argentina
}

\maketitle

\begin{abstract}

We develop a one step matrix method in order to obtain approximate solutions of first order non-linear systems and non-linear ordinary differential equations, reducible to first order systems. We find a sequence of such solutions that converge to the exact solution.   We apply  the method  to different well known examples and check its precision, in terms of local error, comparing it with the error produced by other methods. The advantage of the method over others widely used lies on the great simplicity of its implementation.

\end{abstract}

\section{Introduction}

This paper pretends to be a contribution to methods to find the approximate solutions of nonlinear first order equations (or systems) with given initial values of the form

\begin{equation}\label{1}
{\mathbf y}'(t) = {\mathbf f}({\mathbf y}(t))\,,\qquad {\mathbf y}(t_0)={\mathbf y}_0\,,
\end{equation}
where the prime represents first derivative with respect to the variable $t$. As many higher order ordinary differential equations either linear or non-linear may be written as a system of the form \eqref{1}, our method to obtain  approximate solutions will apply also to these kind of equations.  

Our approach is based in a generalization of the one-step matrix method developed by Demidovich and other authors \cite{DEM,KOT,HSU} some time ago, valid for systems of the form $\mathbf y'(t)=A(t)\mathbf y(t)$. One clear presentation of this method appears in the textbook by Farkas \cite{FAR} and it is interesting to compare it with some related procedures, see for instance \cite{NEP,NEP1}. It is quite important to remark that, while the Demidovich matrix method is applied to linear equations (with variable coefficients), our matrix method is a generalization to non-linear systems.  The advantage that our method may have in comparison of other one step methods lies on its great algorithmic simplicity. In addition, our solutions have a reasonable level of accuracy in few steps and this is workable in a table computer.  This method is quite easily programmable and is very suitable for its use with the package Mathematica. It may be also seen as an alternative to Runge-Kutta and Taylor methods due precisely to its simplicity and precision. The objetive of the present study is to obtain approximate solutions of all kind of systems described by non-linear systems of differential equations (although we may include those linear systems with variable coefficients), including those who appear from physics.

In the derivation of the present approach, our motivation was rooted in the practice of operational numerical calculus.  Nevertheless, we are mainly focused in the mathematical analysis of our method instead of a detailed analysis of the algorithm or CPU times. As for the case of one step Taylor polynomial method, ours shows a great conceptual simplicity. As a consequence, our proposal for obtaining approximate solutions of non-linear systems can be very easily implemented.

The presentation of the method to obtain the approximate solutions is introduced in Section 2. Thus, we have a sequence of approximate solutions, which are defined on a given interval of the real line. This sequence converges uniformly to the exact solution on the given interval as is proven in Section 3, where we also discuss a question of order. It is important to remark that we do not impose any periodicity conditions, so that our results are valid either for periodic or for non-periodic solutions. 

We have applied our method to various examples of two or three dimensional examples.  It is also necessary to test the applicability and accuracy of the method on widely used equations and/or systems. To this end, we have used the van der Pol \cite{VDP},  Duffing \cite{DUF}, Lorentz \cite{LOR} equations, a pseudo diffusive equation depending on a parameter, studied in the standard literature \cite{HAL}, an epidemic equation and the predator-prey Lotka-Volterra \cite{LOT,VOL} equation. We have compared the precision of our solutions with the exact solution, whenever this is known. If not, the comparison is based on Runge-Kutta solutions, for a modern presentation, see \cite{BEL,EGG,KAL}. Also with the widely used Taylor method. 

The analysis of a variety of examples suggest that our method may be more precise for two dimensional systems than for higher dimensional ones, although this is not always exact: we give two examples of three dimensional systems (Lorentz and epidemic equations), for which our method gives different precisions. And it looks like particularly efficient in the case of the pseudo diffusive equation we mentioned earlier, at least when compared with the standard perturbative method studied in \cite{MIC}.

We usually obtain precisions in between of those obtained for the Taylor method of third and fourth order (although the precision depends also on the length of subintervals in which is divided the interval in which we are looking for solutions), which is reasonable when we use a table computer. We close this article with concluding remarks and a conjecture on the Li\'enard equation which is motivated by some of our results and confirmed through numerical experiments.

\section{A matrix method.}

We begin with an initial value problem as given in \eqref{1}. While the function ${\mathbf y}(t)$ is a $\mathbb R^n$ valued function with real variable $t$ of class $C^1$ on the neighbourhood $|t-t_0|\le a$, $\mathbf f(-)$ is a $\mathbb R^n$ valued function with variable in $\mathbb R^n$, which is continuous on $D\equiv \{ \mathbf y \in \mathbb R^n\;\;/\;\; ||\mathbf y - \mathbf y_0||\le d\}$ and  satisfies a Lipschitz condition with respect to $\mathbf y$. Needless to say that $a,d$ are positive constants. 

In relation with the identity \eqref{1}, we shall use either the denomination of ``equation'' or ``system'' indistinctly. In any case, it is well known that the initial value problem \eqref{1} has one unique solution on the interval $|t-t_0|\le \inf(a,d/M)$ with $M:=\sup_D||\mathbf f||$. 

Our objective is to introduce a generalization of a method of solutions of \eqref{1} proposed in \cite{FAR}. This generalization is based in an iterative procedure of numerical integration for equations of the type $\mathbf f(\mathbf y) =A(\mathbf y) \cdot \mathbf y $, where $A(-)$ is an $n\times n$ matrix. With this choice for the function $\mathbf f(-)$, the differential equation in \eqref{1} has the form

\begin{equation}\label{2}
\mathbf y'(t) = A(\mathbf y(t))\cdot \mathbf y(t)\,.
\end{equation}

We assume that the entries of $A(\mathbf y)$ are continuous on $D$. Let us define a uniform partition of the interval $[0,t_N]$ into subintervals $I_k\equiv [t_k,t_{k+1}]$, where $t_k=hk$ with $k=0,1,2,\dots,N$, with $N$ natural and $h=t_N)/N$, $t_N<a$. We have chosen this form of the interval by simplicity, needless to say that if the original interval is somehow else, it always may be transformed into $[0,t_N]$ by a translation. Also, the equal spacing of all subintervals is not strictly necessary, although it simplifies our notation. Conventionally, we may call {\it nodes} to the points $\{t_k\}$. 

We proceed as follows: On each interval $I_k$, we approximate Equation \eqref{2} by

\begin{equation}\label{3}
\mathbf y'_{N,k}(t) = A(\mathbf y^*_{N,k}) \cdot \mathbf y_{N,k}(t)\,.
\end{equation}

At each node, $t_k$, we impose $\mathbf y_{N,k}(t_k) = \mathbf y_{N,k-1}(t_k)$, while $\mathbf y^*_{N,k} \in \mathbb R^n$ is a constant to be determined. This gives the segmentary solution, which has to be of the form $\mathbf y_N(t)\equiv \{\mathbf y_{N,k}(t)\;;\; k=0,1,\dots,N-1\}$. It satisfies

\begin{equation}\label{4}
\mathbf y'_N(t) = A(\mathbf y^*_N) \cdot \mathbf y_N(t)\,,
\end{equation}
where $\mathbf y^*_N$ coincides with $\mathbf y^*_{N,k}$ on each of the $k$-th intervals. 

These segmentary solutions give a sequence of functions $\{\mathbf y_N(t)\}$, $t\in [0,t_N]$, that are approximations to the solution of \eqref{2}. Here, we shall give a method to obtain each of the $\mathbf y_N(t)$ and, in next sections, we shall discuss the properties of the sequence. 

Then, we proceed to an iterative integration of \eqref{4} as follows: Take the first interval $I_0$ and fix an initial value $\mathbf y(t_0)$. This initial value gives the solution $\mathbf y_{N,0}(t)$ on  $I_0$. Thus, we have the value $\mathbf y_{N,0}(t_1)= \mathbf y_{N,1}(t_1)$, which serves as the initial value for the solution on the interval $I_1$. Then, we repeat the procedure in an obvious manner for $I_2$ and so on. For each interval $I_k$ the matrix $A(\mathbf y^*_{N,k})$, which appears in \eqref{3}, is a constant matrix and, therefore, \eqref{3} is a system with constant coefficients. Therefore, the solution of \eqref{3} has the form 

\begin{equation}\label{5}
\mathbf y_{N,k}(t) = \exp \{ A(\mathbf y^*_{N,k})(t-t_k) \} \cdot \mathbf y_{N,k}(t_k)\,,\qquad k=0,1,2,\dots, N-1\,,
\end{equation}
where we have used the notation $\mathbf y_{N,0}(t_0)= \mathbf y(t_0)$. We determine the numbers $\mathbf y^*_{N,k}$ through the following expression:

\begin{equation}\label{6}
\mathbf y^*_{N,k} = \mathbf y_{N,k}\left(t_k+\frac h2 \right) = \exp \left\{ A(\mathbf y_{N,k}(t_k))\,\frac h2  \right\} \cdot \mathbf y_{N,k}(t_k)\,,
\end{equation}
where $k=0,1,2,\dots N-1$. 

Then, the approximate solution $\mathbf y_N(t)$ gives at $t\in I_k$ and on each of the intervals $I_k$:

\begin{equation}\label{7}
\mathbf y_N(t) = \exp \left\{ A(\mathbf y^*_{N,k}) (t-t_k) \right\} \cdot  \prod_{j=0}^{k-1} \exp\left\{ A(\mathbf y^*_{N,j})\,h \right\} \cdot \mathbf y_0\,.
\end{equation}

The determination of the exponential of a matrix may often be rather complicated for large matrices. Then, we may use the Putzer spectral formula. This establish that if $A$ is a constant matrix of order $n\times n$ with eigenvalues $\{\lambda_k\}_{k=1}^n$, then, its exponential verifies the following expression:

\begin{equation}\label{8}
\exp\{ A\,t \} = \sum_{k=1}^n r_k(t)\,P_{k-1}\,,
\end{equation}
with

\begin{equation}\label{9}
P_0\equiv I\,,\qquad P_k=  \prod_{j=1}^k (A-\lambda_j I)\,, \qquad k=1,2,\dots,n-1\,,
\end{equation}
where $I$ is the $n\times n$ identity matrix and the coefficients $r_k(t)$ are to be determined through the following first order system of differential equations:

\begin{eqnarray}\label{10}
r'_1(t) = \lambda_1\,r_1(t)\,, \quad r_1(0)=1\,;  \nonumber \\[2ex] r'_k(t) =\lambda_k\,r_k(t) + r_{k-1}(t)\,, \quad r_k(0)=0\,,
\end{eqnarray}
for $k=2,3,\dots,n$. 

For simplicity, let us consider the particular case, in which $A(\mathbf y^*_{N,j})$ are matrices of order $2\times 2$.  Each of these matrices has two eigenvalues, $\lambda_{1,j}$ and $\lambda_{2,j}$, which may be either different or equal. Let us assume that $\lambda_{1,j} \ne \lambda_{2,j}$. Then,

\begin{eqnarray}\label{11}
\exp\left\{ A(\mathbf y^*_{N,j})\,h \right\} \nonumber\\[2ex]= \frac{1}{\lambda_{1,j} - \lambda_{2,j}} \left\{ (A(\mathbf y^*_{N,j}) - \lambda_{2,j}\,I ) \exp\{\lambda_{1,j}\,h \}   -  (A(\mathbf y^*_{N,j}) - \lambda_{1,j}\,I ) \exp\{\lambda_{2,j}\,h \}  \right\}\,.
\end{eqnarray}

On the other hand, when $\lambda_{1,j}=\lambda_{2,j}=\lambda_j$, we have for the exponential

\begin{equation}\label{12}
\exp\left\{ A(\mathbf y^*_{N,j})\,h \right\} = \exp\{ \lambda_j\,h \}\, \{I + h(A(\mathbf y^*_{N,j}) -\lambda_j\,I)\}\,.
\end{equation}

We conclude here the crude description of the method. In the sequel, we shall show the convergence of the sequence, $\{\mathbf y_N(t)\}$, of approximate segmentary solutions and shall evaluate the precision of the method. 

\section{On the convergence of approximate solutions}

In the previous Section, we have obtained a set of approximate solutions of the initial value problem on a compact interval of the real line. The question is now, assuming we have obtained by the previous method a sequence of solutions. Does this sequence converges to the exact solution in any reasonable sense as the length of the sub intervals, here called $h$, becomes arbitrarily small. To investigate this possibility is the goal of the present Section. Let us go back to \eqref{4} and rewrite it as

\begin{equation}\label{13}
\mathbf y'_N(t)= A(\mathbf y_N) \cdot \mathbf y_N(t) +\eta_N\,,
\end{equation}
so that,

\begin{equation}\label{14}
\eta_N(t) =(A(\mathbf y^*_N) - A(\mathbf y_N)) \cdot \mathbf y_N(t)\,.
\end{equation}

Let us add and subtract $A(\mathbf y^*_N) \cdot \mathbf y^*_N$ in the right hand side of \eqref{14}. Then, let us take the supremum norm on the interval $[t-t_0,t+t_0]$ and use the triangle inequality of the norm, so as to obtain the following inequality:

\begin{equation}\label{15}
||\eta_N|| \le || (A(\mathbf y^*_N) \cdot \mathbf y^*_N -A(\mathbf y_N) \cdot \mathbf y_N)|| + ||A(\mathbf y^*_N)|| \,||\mathbf y^*_N - \mathbf y_N ||\,.
\end{equation}

Then, we apply in \eqref{15} the Lipschitz condition with respect to the variable $\mathbf y$ with constant $K>0$. Then, it comes that

\begin{equation}\label{16}
||\eta_N|| \le (K +|| A(\mathbf y_N^*)||) \, ||\mathbf y_N^*-\mathbf y_N||\,.
\end{equation}

On each one of the intervals $I_k$, let us expand $\mathbf y_N(t)$ in Taylor series around $t_k$. We obtain the following inequality:

\begin{equation}\label{17}
||\mathbf y_N^* - \mathbf y_N(t)|| \le \frac h2 \,||\mathbf y'_n(t_k)|| \le \frac h2\,\max_{t_0\le t\le t_N} ||\mathbf y'_N(t)||\,.
\end{equation}

Since $\mathbf y'_N(t)$ is continuous on the interval $t_0 \le t \le t_N$, the maximum in the right hand side of \eqref{17} exists. Furthermore, $A(\mathbf y)$ is continuous with respect to $\mathbf y$ on the neighborhood $||\mathbf y- \mathbf y_0||\le d$, so that there exists a constant $C>0$ such that $||A(\mathbf y)|| \le C$. Taking norms in \eqref{4}, we have that

\begin{equation}\label{18}
||\mathbf y'_N|| \le ||A(\mathbf y^*_N)||\,||\mathbf y_N|| \le C\, ||\mathbf y_N||\,.
\end{equation}

Equation \eqref{7}  implies that

\begin{equation}\label{19}
\max_{t_0 \le t \le t_N} ||\mathbf y_N(t)|| \le C' \exp\{C(t_N-t_0)\}\,,
\end{equation}
where $C'>0$ is a constant. After (\ref{18}-\ref{19}), we see that $||\mathbf y'(t)||$ is bounded for all $t$ in the interval $[t_0,t_N]$. Consequently after \eqref{16} and (\ref{18}-\ref{19}), we have that

\begin{equation}\label{20}
||\eta_N|| \le \frac h2 \, (K+C)\, C' \, \exp\{C(t_N-t_0)\} =S\,h\,,
\end{equation}
where the meaning of the constant $S>0$ is obvious. 

Next, le us integrate \eqref{13} on the interval $[t_0,t]$. Since for all value of $N$, we use the same initial value $y(t_0)$, we have 

\begin{equation}\label{21}
\mathbf y_N(t) =  \mathbf  y(t_0) + \int_{t_0}^t A(\mathbf y_N(s))\cdot \mathbf y_N(s)\,ds + \int_{t_0}^t \eta_N(s)\,ds\,.
\end{equation}

From \eqref{21}, we have that

\begin{equation}\label{22}
||\mathbf y_{N+M}(t) -\mathbf y_N(t)|| \le \int_{t_0}^t ||A(\mathbf y_{N+M}(s))\cdot \mathbf y_{N+M}(s) - A(\mathbf y_N(s)) \cdot \mathbf y_N(s)||\,ds + \delta_N(t)\,,
\end{equation}
with 

\begin{equation}\label{23}
\delta_N(t) = \int_{t_0}^t ||\eta_{N+M}(s)-\eta_N(s)||\,ds \le 2Sh (t_N-t_0)\,.
\end{equation}

Using the Lipschitz condition in \eqref{22}, we obtain

\begin{equation}\label{24}
||\mathbf y_{N+M}(t) -\mathbf y_N(t)|| \le K \int_{t_0}^t ||\mathbf y_{N+M}(s) -\mathbf y_N(s)|| \, ds + 2Sh (t_N-t_0)\,.
\end{equation}

At this point, we use the Gronwall lema, which states the following

\medskip

{\bf Lemma (Gronwall)}.- Let $f(t): I\longmapsto \mathbb R$ an integrable function on the compact real interval $I$ such that there exists two positive constants $A$ and $B$, with

\begin{equation}\label{25}
0 \le f(t) \le A + B \int_{t_0}^t f(s)\,ds\,, \qquad t_0 \in I\,
\end{equation}
for all $t \in I$. Then,

\begin{equation}\label{26}
f(t) \le A\, e^{B(t-t_0)}\,.
\end{equation}
\hfill $\blacksquare$

\bigskip
Then, we use the Gronwall lema with

\begin{equation}\label{27}
f(t)\equiv ||\mathbf y_{N+M}(t) -\mathbf y_N(t)|| \,, \quad A\equiv Mh := 2S(t_N-t_0)h\,, \quad B\equiv K\,,
\end{equation}
to conclude that

\begin{equation}\label{28}
||\mathbf y_{N+M}(t) -\mathbf y_N(t)|| \le Mh\, e^{K(t-t_0)} \le [M\,e^{K(t_n-t_0)}]\,h =K'\,h\,,
\end{equation}
With $K'>0$ a positive constant. Therefore,

\begin{equation}\label{29}
||\mathbf y_{N+M}(t) -\mathbf y_N(t)|| \longmapsto 0\,,
\end{equation}
as $h\longmapsto 0$. Since the space $C^0[t_0,t_N]$ is complete\footnote{Note that $t_N$ is fixed and $N$ just denotes the number of intervals in the partition or equivalently, the length of $h$.}, \eqref{29} implies the existence of a continuous function $\mathbf z(t): [t_0,t_N] \longmapsto \mathbb R$, such that 

\begin{equation}\label{30}
\mathbf z(t) := \lim_{N\to\infty} \mathbf y_N(t)\,,
\end{equation}
uniformly. 

Now, we claim that $\mathbf z(t)$ is differentiable on $(t_0,t_N)$. Furthermore, $\mathbf z(t)$ is a solution of the differential equation \eqref{2}.

The proof goes as follows: The Lipschitz condition applied to our situation implies that

\begin{equation}\label{31}
||A(\mathbf y_{N+M}(t)) \cdot \mathbf y_{N+M}(t) -A(\mathbf y_N(t)) \cdot \mathbf y_N(t)|| \le K\,||\mathbf y_{N+M}(t)-\mathbf y_N(t)||\,,
\end{equation}
so that $A(\mathbf y_N(t)) \cdot \mathbf y_N(t)$ converges uniformly to $A(\mathbf z(t))\cdot \mathbf z(t)$.  In addition, after \eqref{20}, we have that

\begin{equation}\label{32}
||\eta_N(t)|| \le S h \le S(t_N-t_0)\,.
\end{equation}
Recall that $t_N$ is fixed for whatever value of $N$. Then, taken limits in \eqref{21}, we have

\begin{eqnarray}\label{33}
\mathbf z(t) = \mathbf y(t_0) + \lim_{N\to\infty} \int_{t_0}^t A(\mathbf y_N(s)) \cdot \mathbf y_N(s)\,ds + \lim_{N\to\infty} \int_{t_0}^t \eta_N(s)\,ds \nonumber\\[2ex]   = \mathbf y(t_0) +  \int_{t_0}^t [ \lim_{N\to\infty} A(\mathbf y_N(s)) \cdot \mathbf y_N(s)]\,ds +  \int_{t_0}^t \lim_{N\to\infty}[\eta_N(s)]\,ds\,.
\end{eqnarray}

In the second integral, we may interchange the limit and the integral due to the uniform convergence of the sequence under the integral to its limit. In the case of the second integral, we have used the Lebesgue convergence theorem, which can be applied here due to \eqref{32}.  Since obviously $\lim_{N\to\infty}[\eta_N(s)]=0$, we finally conclude that

\begin{equation}\label{34}
\mathbf z(t)= \mathbf y(t_0) + \int_{t_0}^t A(\mathbf z(s))\cdot \mathbf z(s) \,ds\,.
\end{equation}

From \eqref{34}, we conclude the following:

\medskip

1.- The function $\mathbf z(t)$ is differentiable in the considered interval.

\smallskip
2.- The function $\mathbf z(t)$ is the solution of equation \eqref{2} with initial value $\mathbf z(t_0)=\mathbf y(t_0)$. 

\bigskip

\subsection{A question of order}

The expansion into Taylor series on a neighborhood of $t_k$ of the solutions of equations \eqref{2} and \eqref{3} have these forms, respectively:

\begin{equation}\label{35}
\mathbf y(t_{k+1}) = \mathbf y(t_k) +A(y_k)\,\mathbf y(t_k)\,h +\frac 12 (A^2(\mathbf y(t_k))\cdot \mathbf y(t_k))h^2  + \frac 12 \frac{d}{dt} [A(\mathbf y(t_k)) \cdot \mathbf y(t_k)] h^2 + O(h^3)\,.
\end{equation}

Taking into account that

\begin{equation}\label{36}
\frac{d}{dt}\,A(\mathbf y(t)) = \sum_{j=1}^n \frac{\partial}{\partial\,y_j}\, A(\mathbf y(t)) \,\frac{d}{dt}\, y_j(t)\,,
\end{equation}
where $y_j$ is the $j$-th component of $\mathbf y$, equation \eqref{35} becomes:

\begin{eqnarray}\label{37}
\mathbf y(t_{k+1}) = \mathbf y(t_k) +A(y_k)\,\mathbf y(t_k)\,h +\frac 12 (A^2(\mathbf y(t_k))\cdot \mathbf y(t_k))h^2 \nonumber\\[2ex] + \frac12 \left( \sum_{j=1}^n \frac{\partial}{\partial\,y_j}\, A(\mathbf y(t_k)) \,\frac{d}{dt}\, y_j(t_k) \right) \cdot \mathbf y(t_k) \,h^2 + o(h^3)\,.
\end{eqnarray}

On the $k$-th interval, equation \eqref{4} takes te formula

\begin{equation}\label{38}
\mathbf y_N(t_{k+1}) = \mathbf y_N(t_k) + A(\mathbf y^*_{N_k}) \cdot \mathbf y_N(t_k)\,h +\frac 12\, A^2(\mathbf y^*_{N_k}) \cdot \mathbf y_N(t_k)\,h^2 +o(h^3)\,.
\end{equation}

Then, we may proceed to expand into Taylor series the matrix $A(\mathbf y^*_{N_k})$ on a neighborhood of $\mathbf y_{N_k}$. Taking into account \eqref{5} and after some simple calculations, we obtain:

\begin{equation}\label{39}
A(\mathbf y^*_{N_k}) = A(\mathbf y_N(t_k) ) + \sum_{j=1}^n \frac{\partial}{\partial \,y_j}\, A(\mathbf y_N(t_k) ) \,(y_{N,j}(t_k+h/2)- y_{N,j}(t_k))\,,
\end{equation}
where $y_{N,j}(t)$ is the $j$-th component of $\mathbf y_N(t)$. A first order expansion on the last factor on the right hand side of \eqref{39} gives

\begin{equation}\label{40}
y_{N,j}(t_k+h/2)- y_{N,j}(t_k) = \frac 12\,\frac{d}{dt}\, y_{N,j}(t_k) \, h + o(h^2)\,,
\end{equation}
so that using \eqref{40} in \eqref{39}, we have

\begin{equation}\label{41}
A(\mathbf y^*_{N_k}) = A(\mathbf y_N(t_k) ) + \frac h2 \sum_{j=1}^n \frac{\partial}{\partial \,y_j}\, A(\mathbf y_N(t_k) ) \, \frac{d}{dt} \, y_{N,j}(t_k)\,.
\end{equation}

Then, we replace \eqref{41} into \eqref{38}, and performing some simple manipulations, taking into account that up to second order in $h$:

\begin{equation}\label{42}
\frac12 \, A^2(\mathbf y^*_{N_k}) \cdot \mathbf y_N(t_k)\,h^2 \approx \frac 12\,A^2(\mathbf y(t_k)) \cdot \mathbf  y_N(t_k)\, h^2\,,
\end{equation}
we finally obtain that

\begin{eqnarray}\label{43}
\mathbf y(t_{k+1})= \mathbf y_N(t_k) + A(\mathbf y_N(t_k)) \cdot \mathbf y_N(t_k)\,h +\frac{h^2}{2}\sum_{j=1}^n \,\frac{\partial}{\partial\,y_j} A(\mathbf y_n(t_k))\,\frac{d}{dt}\, y_{N,j}(t_k) \nonumber \\[2ex] + \frac12 \,A^2(\mathbf y_N(t_k)) \cdot \mathbf y_N(t_k)\,h^2 +o(h^3)\,.
\end{eqnarray}

The advantage that \eqref{43} offers with respect to \eqref{37} is that in \eqref{43} the terms up to second order in $h$ are correctly shown.

\subsubsection{Going beyond second order}

In (3.1), we have found the solution up to second order in $h$. Would we obtain a better precision for third or higher order keeping at the same time the simplicity of the method? First of all our construction is based on equation \eqref{3}, which is not longer valid if we require an approximation of order higher than two.  To fix ideas, let us take  ${\mathbf y}(t)=(y_1(t),y_2(t))$ bidimensional for simplicity, a choice which does not affect to our argument. Let us expand $A({\mathbf y}(t))$ around ${\mathbf y}^*(t)=(y_1^*,y_2^*)$ (where we have omitted the sub-indices $N,k$ for simplicity) and take one more term in the Taylor span. The result is

\begin{equation}\label{44}
A({\mathbf y}(t)) = A({\mathbf y}^*(t)) + \frac{\partial}{\partial y_1}\, A({\mathbf y}^*(t)) (y_1(t)- y_1^*) + \frac{\partial}{\partial y_2}\, A({\mathbf y}^*(t)) (y_2(t)- y_2^*)\,.
\end{equation}

Then, instead of the approximation \eqref{3}, we have the following, where again we have suppressed the indices $N$ and $k$:

\begin{equation}\label{45}
{\mathbf y}'(t) = \left[ A({\mathbf y}^*(t)) + \frac{\partial}{\partial y_1}\, A({\mathbf y}^*(t)) (y_1(t)- y_1^*) + \frac{\partial}{\partial y_2}\, A({\mathbf y}^*(t)) (y_2(t)- y_2^*) \right] \cdot {\mathbf y}(t)\,.
\end{equation}

The solution to be determined is just ${\mathbf y}(t)=(y_1(t),y_2(t))$, which now is a part of the Ansatz \eqref{45}. One possibility to solve this contradiction is to proceed with the following span on each of the $I_k$ intervals:

\begin{equation}\label{46}
y_1(t)= y_1(t_k) + \frac{\partial y_1}{\partial t}(t_k) \,(t-t_k) +\dots\,,
\end{equation}
and same for $y_2(t)$. If we use these manipulations in \eqref{3}, we obtain an equation of the type:

\begin{equation}\label{47}
{\mathbf y'_{N,k}}(t) = G_{N,k}(t) \cdot {\mathbf y_{N,k}}(t) \,,
\end{equation}
with 

\begin{equation}\label{48}
G_{N,k}(t) = A({\mathbf y^*_{N,k}}(t)) + \frac{\partial }{\partial y_1}\, A({\mathbf y^*_{N,k}}(t)) (y_1(t)-y_1^*) + \frac{\partial }{\partial y_2}\, A({\mathbf y^*_{N,k}}(t)) (y_2(t)-y_2^*) + \dots\,.
\end{equation}

Obviously, system \eqref{48} is non-autonomic. We have to use a new approximation of $G_{N,k}(t) $ by a constant on the interval $I_k$ and re-start again. 

As we see, advancing to just one higher order of accuracy destroys the simplicity of the method which is one of its more interesting added values. Therefore, we cannot consider going to higher orders an advantage. It may slightly improve the precision at the price of destroying the efficiency and simplicity of the method.

\subsection{Some examples}

\begin{itemize}

\item{{\bf The van der Pol equation}

The van der Pol equation

\begin{equation}\label{49}
y''(t) +\mu(y^2(t)-1)y'(t) + y(t)=0
\end{equation}
is a particular case of the Li\'enard equation, 

\begin{equation}\label{50}
y''(t)+f(y)\,y'(t) +g(y)=0\,.
\end{equation}
which will be discussed in the Appendix. In the van der Pol equation, we have obviously that $f(y)=\mu(y^2(t)-1)$ and $g(y)\equiv y$. This equation can be easily written in the matrix form \eqref{2}  by writing $y_1(t)\equiv y(t)$ and $y_2(t)\equiv y'(t)$, as

\begin{equation}\label{51}
\left(\begin{array}{c} y'_1(t)\\[2ex] y'_2(t) \end{array}  \right) =  \left( \begin{array}{cc} 0 & 1 \\[2ex] -1 & \mu(1-y_1^2) \end{array} \right)   \left(\begin{array}{c} y_1(t)\\[2ex] y_2(t) \end{array}  \right)\,.
\end{equation}

Our goal is to compare the precision of our method with a reference solution. No explicit solutions to the van der Pol equation \eqref{49} are known, so that we use the Runge-Kutta solution of eight order, $y_{rk}(t)$, as reference solution (alternatively, one may consider a Taylor solution of eighth order, which has a comparable precision). We compare the precision of our method with the precision of the solutions as obtained by a third or fourth  order Taylor method. As a measure of the error, we use 

\begin{equation}\label{52}
e_h:= \frac 1N \sum_{j=0}^{N-1} (y_{rk}(t_j)-y(t_j))^2\,,
\end{equation}
where $y(t)$ is the solution obtained using our method or any other, like the Taylor method. Then, we need to use given values of the parameters and inicial conditions. In  Table 1 below, we compre the errors produced by our method as compared with the error given with the use of third and fourth order Taylor method for different values of the interval width $h$ for $T=20$, $\mu=1/2$ and the initial conditions $y(0)=0$ and $y'(0)=2$. 

\vskip1cm

\centerline{$
\begin{array}
[c]{cccc}
h & {\rm Matrix} & {\rm Taylor\, 3^{\rm rd}}  & {\rm Taylor\, 4^{\rm rd}}\\[2ex]
10^{-4} & 2.63\, 10^{-11} & 8.60\,10^{-7} & 6.40 \,10^{-7} \\
10^{-3} & 2.08\, 10^{-11}  & 8.30\, 10^{-9} & 6.45\, 10^{-9} \\
10^{-2}  & 1.46\, 10^{-8} & 9.11\, 10^{-9} & 6.60\, 10^{-11} \\
10^{-1} & 2.04\, 10^{-5} & 1.17\, 10^{-4} & 1.67\, 10^{-7} \\
2\, 10^{-1} & 2.30\, 10^{-4} & 2.23\, 10^{-3} & 7.94\, 10^{-6} \\
5\, 10^{-1} & 8.10\, 10^{-3} & 1.49\, 10^{-1} & 2.75\, 10^{-3}
\end{array}
$}

\medskip
TABLE 1.- Values for the error $e_h$ for our matrix method and the Taylor method of orders three and fourth for distinct values of $h$ for the van der Pol equation.

\vskip 1cm

It is clear that our matrix method has a precision in between those of the third and fourth orden Taylor method. These results are just an example of the results obtain in multiple numerical experiments, we have performed showing essentially the same result. However, Table 1 as well as other numerical experiments show and important tendency: the lower $h$ the better our precision as compared with the precision given by the Taylor method. The reason is clear, the smaller $h$ the bigger the number of operations needed to obtain the approximate solution. Our matricidal method requieres less operations than the Taylor method, so that our precision gets better as $h$ gets smaller. 

In Appendix B, we give our source code when our method is to be applied in the present case.

}

\item{{\bf The Duffing equation}

This is another second order non linear equation, which has the following form:

\begin{equation}\label{53}
y''(t) + y(t) + y^3(t)=0\,,
\end{equation}
where we have omitted the term in the first derivative of the indeterminate, $y'(t)$. For the Duffing equation, there are explicit solutions in terms of the Jacobi elliptic functions. For instance, being given the initial conditions $y(0)=1$ and $y;(0)=0$, we have the following solution:

\begin{equation}\label{54}
y_e(t) = -i \sqrt{1+k} \,{\rm sn}\,(u;m)\,,
\end{equation}
where ${\rm sn}\,(u;m)$ denotes the elliptic sine. The arguments in \eqref{49} denote the following:

\begin{eqnarray}\label{55}
u= \frac 1{\sqrt 2} \, (t^2(1-k) + 2t(c_2-kc_1) + (1-k)c_2^2)\,, \nonumber\\[2ex]  m= \frac{1+k}{1-k}\,,\qquad k=\sqrt{1+2c_1}\,. 
\end{eqnarray}

The value of the constants included in \eqref{55} are $c_1= 1.5$ and $c_2= 1.1920055$. The Duffing equation may be written in matrix form, if we write again $y_1(t)\equiv y(t)$ and $y_2(t)\equiv y'(t)$, as

\begin{equation}\label{56}
\left(\begin{array}{c} y'_1(t)\\[2ex] y'_2(t) \end{array}  \right) =  \left( \begin{array}{cc} 0 & 1 \\[2ex] -(1+y_1^2) & 0 \end{array} \right)   \left(\begin{array}{c} y_1(t)\\[2ex] y_2(t) \end{array}  \right)\,.
\end{equation}

We define the error $e_h$ as in \eqref{52}, where we replace $y_{rk}(t)$ by the exact solution, $y_e(t)$, which does exist in the present case. Also $y(t)$ is the solution for which we want to compare its precision with the exact solution, in our case the solutions obtained by our matrix method as well as the third or fourth order Taylor solutions. The errors produced by each method are given on Table 2 below. 

\vskip1cm

\centerline{$
\begin{array}
[c]{cccc}
h & {\rm Matrix} & {\rm Taylor\, 3^{\rm rd}}  & {\rm Taylor\, 4^{\rm rd}}\\[2ex]
10^{-4} & 1.69\, 10^{-9} & 5.60\,10^{-5} & 4.13 \,10^{-5} \\
10^{-3} & 1.03\, 10^{-10}  & 5.06\, 10^{-7} & 4.24\, 10^{-7} \\
10^{-2}  & 9.61\, 10^{-9} & 2.47\, 10^{-8} & 4.24\, 10^{-9} \\
10^{-1} & 1.16\, 10^{-5} & 4.92\, 10^{-4} & 4.12\, 10^{-8} \\
2\, 10^{-1} & 1.18\, 10^{-4} & 6.78\, 10^{-3} & 7.56\, 10^{-6} \\
5\, 10^{-1} & 82.00\, 10^{-3} & 1.04\, 10^{-1} & 7.73\, 10^{-3}
\end{array}
$}

\medskip
\medskip
TABLE 2.- Values for the error $e_h$ for our matrix method and the Taylor method of orders three and fourth for distinct values of $h$ for the Duffing equation.

\vskip1cm

We see that the results are quite similar to those obtained with the van der Pol equation. Similarly, we have made some numerical experiments that confirm these results. 

}

\item{ {\bf The Lorenz equation}

The Lorenz equation, which is a model for the study of chaotic systems has been introduced in the study  of atmospheric behaviour \cite{LOR}. This equations arises in many problems of physics, where chaoticity is present \cite{HAK,GOR,COU,MIS}. The Lorenz equation is usually written in matrix form and is three-dimensional:

\begin{equation}\label{57}
\left(\begin{array}{c} y'_1(t) \\[2ex] y'_2(t) \\[2ex] y'_3 (t) \end{array}  \right) = \left(\begin{array}{ccc} -a & a & 0 \\[2ex] b-y_3(t) & -1 & 0 \\[2ex] y_2(t) & 0 & -c \end{array} \right)  \left(\begin{array}{c} y_1(t) \\[2ex] y_2(t) \\[2ex] y_3 (t) \end{array}  \right)\,,
\end{equation}
$a$, $b$ and $c$ being positive constants. 

This system is very sensitive to  the particular choice of the parameters and of initial values, as may numerical experiments show. Based in these experiments, we have choosed the following values for the parameters, $a=10$, $b=9.996$ and $c=8/3$, the fixed point $(y_1,y_2,y_3) =(4.88808,4.88808,8.996)$ is an attractor. We proceed as in the previous case and compare solutions of the errors \eqref{52} resulting of the use of the Taylor method of order two and this Matrix method using as reference the solution obtained with the Ruge Kutta method of eighth order. We use $T=10$ and the solution with initial values $(y_1(0),y_2(0),y_3(0))= (1,0,1)$. We have obtained a table (Table 3) of errors given in terms of the interval width $h$ as

\vskip1cm

\centerline{$
\begin{array}
[c]{ccccc}
h & {\rm Matrix} & \text{Taylor second order} & \text{Taylor third order} & \text{Taylor fourth order}\\[2ex]
10^{-2} & 5.2\,10^{-6} & 6.73\,10^{-5} & 3.6\,10^{-8}  &  3.1\, 10^{-10} \\
10^{-1} & 5.3\,10^{-2}  & 7.1 \, 10^{-1} & 1.8\, 10^{-1} & 5.2\, 10^{-2} \\
2\, 10^{-1}  & 3.7 \,10^{-1}  & {\rm error}  &  {\rm error}  & {\rm error}
\end{array}
$}

\medskip
\centerline{TABLE 3.- Values of the precision $e_h$ in terms of $h$ for $T=10$, $y_1(0)=1$, $y_2=0$ and $y_3(0)=1$.}

\vskip1cm

The word ``error'' written in two entries in Table 2 means that the error that appears if the interval width is of the order of 0.2, when we use the Taylor method, is incontrollable. We also see that the precision obtained with the matrix method clearly improves the precision by the Taylor method the larger the length of the subintervals.

\item{{\bf Neutral damping equation}.

This equation has been discussed in the literature \cite{HAL,MIC} and is

\begin{equation}\label{58}
 x''(t) +\varepsilon\, ( x'(t))^2 +x(t)=0\,,
\end{equation}
where the tilde means derivative with respect to the variable $t$ and $\varepsilon$ is a real parameter. As in previous cases, let us define $y(t):=  x'(t)$, so that \eqref{58} can be written in the following form:

\begin{equation}\label{59}
(-x+\varepsilon y^2)\,dx -y\,dy=0\,.
\end{equation}

This equation is integrable with integrating factor:

\begin{equation}\label{60}
\mu(x,y)\equiv e^{2\varepsilon x}\,.
\end{equation}

It is readily shown that equation \eqref{60} admits the following first integral:

\begin{equation}\label{61}
f(x,y) =  \left[\frac 12 \,y^2 + \frac{1}{4\varepsilon^2}\,(2\varepsilon x -1) \right] e^{2\varepsilon x}
\end{equation}

If we write \eqref{58} in  the standard matrix form as

\begin{equation}\label{62}
\left(\begin{array}{c} x' \\[2ex]  y' \end{array} \right) = \left( \begin{array}{cc} 0 & 1 \\[2ex] -1 & -\varepsilon y \end{array} \right)  \left( \begin{array}{c} x\\[2ex] y \end{array} \right)\,,
\end{equation}
we readily observe that the only fixed point is the origin. It is also the origin the point at which the minimum of the first integral \eqref{61} lies. All orbits are periodic around the origin. 

In Figure 1, we have depicted the periodic solutions in phase space. Note that, although equation \eqref{58} may look like dissipative, it is not as Figure 1 manifests. 

Now, let us repeat the comparison between errors given by our matrix method with the errors given by the Taylor method at some low orders. As always, we use \eqref{47} as the definition of the error and the numerical Runge-Kutta solution of eighth order as the reference solution. We obtained the following results, which appear in Table 4, where $h$ is, as always, the distance between two consecutive nodes:

\vskip1cm

\centerline{$
\begin{array}[c]{cccc}
h & {\rm Matrix} & {\rm Taylor\, 3^{\rm rd}}  & {\rm Taylor\, 4^{\rm rd}}\\[2ex]
10^{-4} & 2.2\, 10^{-14} & 4.2\, 10^{-14}  & 3.2 \,10^{-11} \\
10^{-3} & 1.9\,10^{-14} & 4.1\,10^{-11} &  3.1\, 10^{-11} \\
10^{-2} & 9.0\, 10^{-14} & 4.6\, 10^{-13} & 4.0\,10^{-12} \\
10^{-1} & 1.6\, 10^{-9} & 2.7\, 10^{-7}  & 8.2\, 10^{-11} \\
2.10^{-1} & 2.6\, 10^{-8} & 8.6\, 10^{-6} & 1.0\, 10^{-9} \\
5.10^{-1} & 1.1 \, 10^{-6} & 7.5\, 10^{-4} & 6.1\, 10^{-6}
\end{array}
$}
\medskip
TABLE 4. Values of the error in terms of $h$ for $\varepsilon=0.1$, $T=2\pi$ (see \eqref{47}) and the initial values $x(0)=1$, $ x'(0)=0$. 

\vskip1cm 

We observe that the precision of the matrix method is much higher than the precision of the third and forth order Taylor method for a distance between nodes $h\le 0.01$. Moreover, we have to underline that our matrix one step method is much simpler to programming that the Taylor method as the reader can easily convince him/herself using any of these examples.  

There is another method based in the theory of perturbations in order to obtain approximate solutions to \eqref{53}, which receives the name of Lindstedt-Poincar\'e. It is described in \cite{MIC}. It consists in a series in terms of $\varepsilon$ of the form:

\begin{equation}\label{63}
x(t) = x_0(t) +\varepsilon\, x_1(t) + \varepsilon^2\,x_2(t) +\dots\,.
\end{equation}

Coefficients $x_i(t)$ may be obtain iteratively, once we have fixed the initial values. For instance, for $x(0)=1$ and 
$ x'(0)=0$, we obtain:

\begin{eqnarray}\label{64}
x_0(t)= \cos\omega t\,,\qquad x_1(t)=\frac 16 \left(-3+4 \cos \omega t- \cos 2\omega t \right)\,,\nonumber\\[2ex]
x_2(t)= \frac 13 \left(-2 +\frac{61}{24}\, \cos\omega t -\frac 23 \,\cos 2\omega t + \frac 18\, \cos 3\omega t \right)\,, \nonumber\\[2ex]
\omega= 1-\frac 16\,\varepsilon + o(\varepsilon^3)\,.
\end{eqnarray}

We may also evaluate the error for this perturbative method, which is independent of any division of the interval, in which the solution is considered, into subintervals. This error is $1.1\,10^{-4}$, which is obviously higher to those obtained using any of the numerical method considered. 

We finish the present example by proposing another approach to an approximate solution. By either method, matrix or Taylor, we obtain on each of the nodes $\{t_k\}$ a value, $x_k$, of the approximate solution. Let us interpolate each interval by means of cubic splines, so that we obtain a segmentary approximation by cubic polynomials \footnote{In \cite{LG}, we have already proposed segmentary cubic solutions and studied their properties, although in \cite{LG} they were not necessarily cubic splines.}. Let us assume that the cubit interpolating solution is $s(t)$, then the error of this solution on an interval $T$, with respect to the exact solution is given by (recall that cubic splines admit first and second continuous derivatives at the nodes)

\begin{equation}\label{65}
e=\frac 1T \int_0^T (s''(t)+\varepsilon\,[s'(t)]^2+s(t))^2\,dt\,.
\end{equation} 

The form of the error for the Taylor method is also given by \eqref{60}. The resulting errors appear in the following table (Table 5):

\vskip1cm 

\centerline{$
\begin{array}[c]{ccccc}
h & {\rm Spline} & {\rm Taylor\, 2^{\rm rd}} & {\rm Taylor\, 3^{\rm rd}}  & {\rm Taylor\, 4^{\rm rd}}\\[2ex]
10^{-4} & 7.2\, 10^{-16} & 4.7 \,10^{-16} & 6.7\,10^{-16} & 6.5\,10^{-16}\\
10^{-3} & 1.8\,10^{-14} & 1.8\,10^{-14} & 1.8\,10^{-14} & 1.8\,10^{-14}\\
10^{-2} & 1.9\,10^{-10} & 1.9\,10^{-11} & 1.1\,10^{-10} & 1.1\,10^{-10} \\
10^{-1} & 4.6\,10^{-6} & 6.4\,10^{-6} & 4.6\,10^{-6} & 4.6\,10^{-6} \\
2.10^{-1} & 6.1\,10^{-5} & 9.8\,10^{-5} & 6.0\,10^{-5} & 6.0\,10^{-5}\\
5.10^{-1} & 4.8\,10^{-3} & 1.2\,10^{-3} & 4.6\,10^{-3} \ & 4.9\,10^{-3}
\end{array}
$}

\medskip

TABLE 5.- Comparison between errors by the cubic spline and Taylor methods. 

\vskip1cm

Finally, for the perturbative method the error obtained is $7.4\,10^{-5}$. Concerning the conservation of the constant of motion, we measure its dispersion by means of the following parameter, $e_f$, defined as

\begin{equation}\label{66}
e_f:= \frac 1T \int_0^T(f(x_0,y_0)-f(x,y))^2\, dt\,,
\end{equation}
where $f(x,y)$ should be calculated using different approximations such as Taylor, matrix and the analytic approximate solution as obtained by the perturbative Lindstedt-Poincar\'e method mentioned earlier. The point $(x_0,y_0)$ gives the chosen initial conditions. In the latter case, we have obtained $e_f = 1.1\,10^{-4}$, in all others, we always got $e_f<10^{-8}$.

\begin{figure}
\centering
\includegraphics[width=0.5\textwidth]{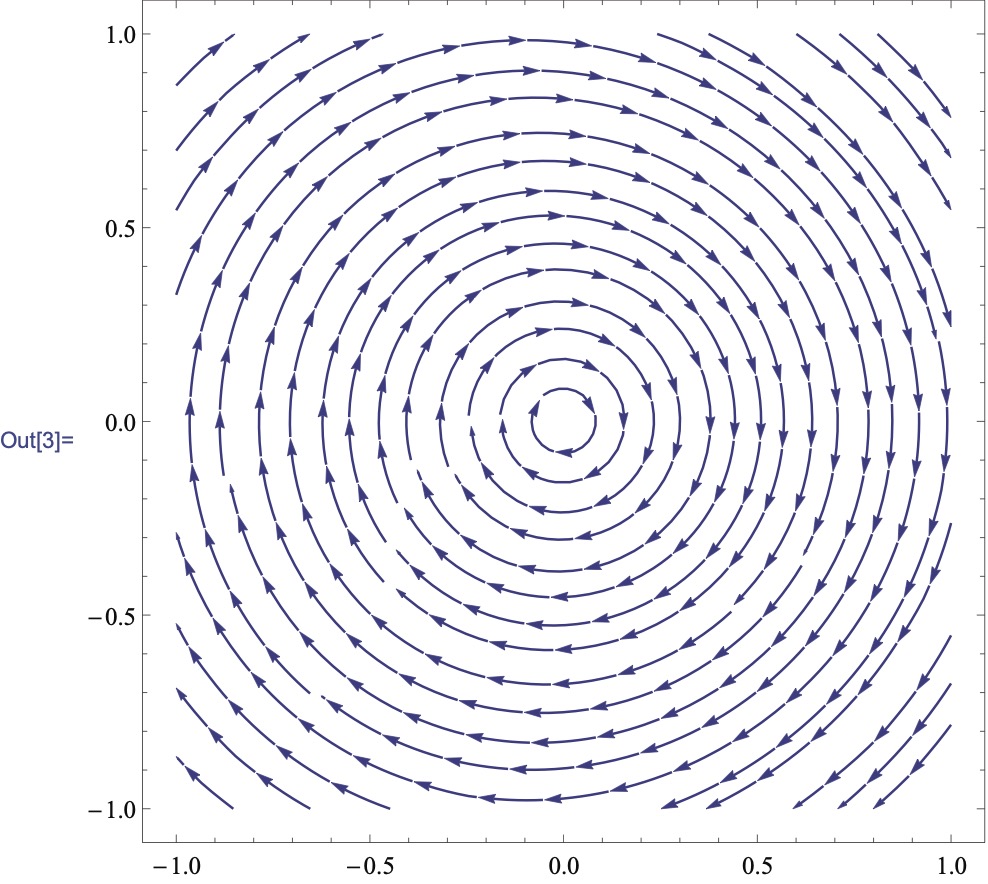}
\caption{\small Periodic orbits around the origin for the neutral damping equation. The horizontal coordinate represents the values of $x$, while vertical coordinate gives the values of $y=x'$. Note that these periodic orbits are represented in phase space.
\label{Figure1}}
\end{figure}
}

\item{{\bf An epidemic equation}

A model for an epidemia has been proposed as early as in 1927 \cite{KK,STR,TC}. If $x(t)$, $y(t)$ and $z(t)$ are the number, at certain time $t$, of healthy, sick and dead persons, respectively, in some society, the model assumes that these functions satisfy the following non-linear system:

\begin{eqnarray}\label{67}
\dot x(t) &=& -a x(t)\,y(t)\,, \nonumber\\[2ex] \dot y(t) &=& a x(t)\,y(t) -b y(t)\,, \nonumber\\[2ex] 
\dot z(t) &=& by(t)\,,
\end{eqnarray}
where $a$ and $b$ are positive constants and the dot means derivation with respect to time $t$. The model assumes infection of healthy persons from sick persons. The latter died after some time. Note that, in the studied time interval, the sole cause of population dynamics is the epidemic. For obvious reasons, we consider only positive solutions. The fixed points for \eqref{67} have the form $(\alpha,0,\beta)$ with $\alpha,\beta>0$. 

Let us consider the following vector field, also called the {\it flux} of \eqref{62}, 

\begin{equation}\label{68}
F(t):=(-a x(t)\,y(t),a x(t)\,y(t) -b y(t),by(t)) \,.
\end{equation}

Note that $F_x<0$. For $x(t)<b/a$, we have $F_y<0$, while $x>b/a$, it results that $F_y>0$. Thus, the fixed point is inestable if $\alpha >b/a$. This is a necessary condition for the existence of the epidemic. If $\alpha <b/a$, then the fixed points are asymptotically stable as depicted in Figure 2, where we have chosen $a=0.0005$ and $b=0.1$. Different curves in Figure 2 represent different initial conditions.

\begin{figure}
\centering
\includegraphics[width=0.5\textwidth]{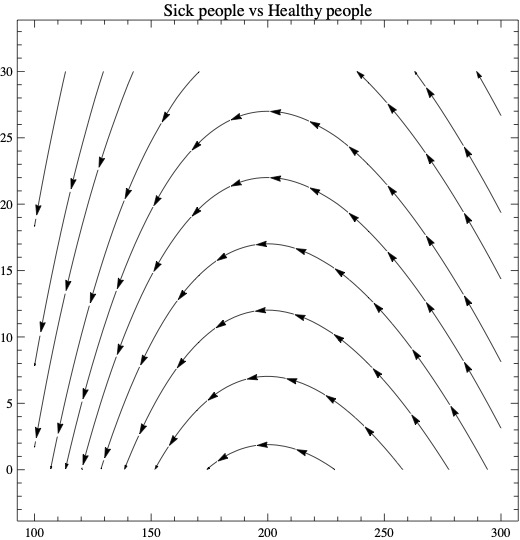}
\caption{\small Asymptotic stability of fixed points with $a=0.0005$ and $b=0.1$ in \eqref{67}. The horizontal line represent the scaled number of healthy people, while the vertical figures give the scaled number of sick people. The arrow means the direction of time. Different curves are obtained using different initial conditions. Observe the presence of a maximum at a point which is independent on the initial conditions.
\label{Figure2}}
\end{figure}

Let us write \eqref{67} in matricial form as

\begin{equation}\label{69}
\dot X(t) =A(x,y,z)\cdot X(t)\,,
\end{equation}
where,

\begin{equation}\label{70}
X(t):= \left(\begin{array}{c} x(t)\\[2ex] y(t)\\[2ex] z(t) \end{array}  \right)\,, \qquad A(x,y,z):= \left( \begin{array}{ccc}  0 & -ax & 0 \\[2ex] 0 & ax-b & 0 \\[2ex] 0 & b & 0 \end{array} \right)\,.
\end{equation}

In Figure 3, we obtain the curve giving the total number of infected people with time. After a maximum, the number of sick people decays quickly. We have used the values for the parameters $a=0.0005$ and $b=0.1$ and the initial values 
$x(0)=300 >b/a$, $y(0)=20$, $z(0)=0$. 

\begin{figure}
\centering
\includegraphics[width=0.5\textwidth]{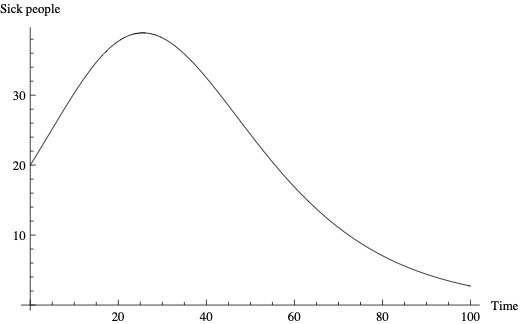}
\caption{\small Number of infected people in relation with time. Observe the existence of a maximum. After the maximum the curve decreases steeply. 
\label{Figure3}}
\end{figure}

The form of the error for the method is obtained using \eqref{65} again. Next in Table 6, we compare the errors for given values of $h$ of our Matrix Method as compared to second and third order Taylor. 

\vskip1cm 

\centerline{$
\begin{array}[c]{cccc}
h & \text{Matrix Method} & {\rm Taylor\, 2^{\rm rd}} & {\rm Taylor\, 3^{\rm rd}}  \\[2ex]
0.01 & 2.2\, 10^{-11} & 4.0 \, 10^{-11} & 1.9\, 10^{-11} \\
0.1 & 2.0\, 10^{-8} & 7.3 \, 10^{-8} & 2.0 \, 10^{-11} \\ 
1.0 & 2.1\, 10^{-4} & 7.2 \, 10^{-4} & 1.7\,10^{-7} \\
2.0 & 3.6\, 10^{-3} & 1.1\, 10^{-2} & 1.2\, 10^{-5}
\end{array}
$}

\medskip

TABLE 6.- Comparison between the errors for the Matrix Method and the Taylor method in the case of the epidemic model. 	

\vskip1cm

We see that the level of error is similar, although our method keeps the advantage  of needing much less arithmetics than any Taylor method. 
}

\item{{\bf Lotka-Volterra equation.}

Models in population dynamics, as for instance the predator-prey competition, were independently developed by the american biologist A.J. Lotka \cite{LOT} and the Italian mathematician V. Volterra \cite{VOL}, see modern references for the Lotka-Volterra equation in \cite{FAR,CI,MUR,LV}. The most general form of the Lotka-Volterra equation has the form

\begin{eqnarray}\label{71}
\dot x_1= x_1(\varepsilon_1 -a_{11}\,x_1 - a_{12}\, x_2)\,, \nonumber\\[2ex]
\dot x_2 =x_2(\varepsilon_2 -a_{12}\,x_1 - a_{22}\, x_2)\,,
\end{eqnarray}
where $x_1$ and $x_2$ are functions of time $t$ and $\varepsilon_i$, $a_{ij}$, $i,j=1,2$ are constants. Following \cite{CI}, we consider here a simpler version, yet non-linear, of \eqref{71}, which is

\begin{eqnarray}\label{72}
\dot x= a\,x-b\,xy\,,\nonumber\\[2ex] \dot y = d\,xy - c\,y\,,
\end{eqnarray}
and the initial conditions $x(t_0)=x_0$, $y(t_0)= y_0$.

We may consider the solutions, $x=x(t)$, $y=y(t)$ of \eqref{72} as the equations determining a parametric curve on the $x-y$ plane. Then, by elimination of $t$ and integration, we obtain:

\begin{equation}\label{73}
C(x,y)= x^cy^a \,\exp\{-(b+d)x\}\,.
\end{equation}

Note that $C(x,y)$ is a constant on each curve solution (constant of motion) and its value over each solution is determined via the initial conditions. Equations \eqref{71} have two fixed points, which are $(0,0)$ and $(c/d,a/b)$. For $x>0$ and $y>0$ all solutions are periodic. 

To apply the method proposed in the present article to this system, let us write it in matrix form as

\begin{equation}\label{74}
\dot X =A(x,y)\,X\,,
\end{equation}
where,

\begin{equation}\label{75}
X(t)\equiv \left(\begin{array}{c} x(t)\\[2ex] y(t) \end{array} \right)\,,\qquad A\equiv   \left(\begin{array}{cc} a-by & 0 \\[2ex] 0 & dx-c \end{array} \right)\,.
\end{equation}

In order to estimate the precision of the method, we need to choose values for the initial values as well as for the parameters $a$, $b$, $c$ and $d$. We have use several choices and obtain in all of them similar values. to show a table comparing the precision of our method with some others, let us choose as the values of the parameters, $a=1.2$, $b=0.6$, $c=0.8$and $d=0.3$. As the initial values, let us choose, $x_0=y_0=3$. In addition, we have to take an integration time, in our case we took $T=20$, which approximately accounts for three periods. 

In order to determine the error produced by our matrix method is determined by the formula \eqref{65} above. We compare this error with those of second and third order Taylor method as compared to the numerical solution obtained by a forth order Runge-Kutta. These errors are shown in Table 7.

\vskip1cm

\centerline{$
\begin{array}[c]{cccc}
h & \text{Matrix Method} & {\rm Taylor\, 2^{\rm rd}} & {\rm Taylor\, 3^{\rm rd}}  \\[2ex]
0.01 & 4.1\, 10^{-8} & 3.5 \, 10^{-8} & 8.2\, 10^{-13} \\
0.1 & 3.5\, 10^{-8} & 5.4 \, 10^{-8} & 2.3 \, 10^{-9} \\ 
0.2 & 2.6\, 10^{-3} & 5.4 \, 10^{-3} & 1.0\,10^{-5} 
\end{array}
$}

\medskip

TABLE 7.- Comparison between the errors for the Matrix Method and the Taylor method in the case of the Lotka-Volterra equation. 	

\vskip1cm

We see that the precision of our method is equivalent to those of a second order Taylor, with much less arithmetic operations. 

}

\item{{\bf On the possibility of extending the method to PDE: The Burger's equation.}

Can we extend this precedent discussion to partial differential equations admitting separation of variables? One possible example would have been the convection diffusion equation in one spatial dimension \cite{HWX}:

\begin{equation}\label{76}
\frac{\partial}{\partial t} u(x,t)= \frac{\partial}{\partial x} \left[ D(u)\,\frac{\partial}{\partial x}\, u(x,t) \right] +Q(x,u)\, \frac{\partial}{\partial x}\,u(x,t) + P(x,u)\,.
\end{equation}

Separation of variables for \eqref{76} is discussed in \cite{HWX}. In general, our method is not applicable here since the resulting equations after separation of variables are not of the form \eqref{1}. Nevertheless, another point of view is possible. Assume we want to obtain approximate solutions of an equation of the type \eqref{76} on the interval $[0,X]$ under the conditions $u(0,t)=0$, $u(X,t)=0$, $u(x,0)=h(x)$, $h(x)$ being a given smooth function and $u(0,0)=u(X,0)=0$. On the interval $[0,X]$, we define a uniform partition of width $h:=X/n$ and nodes $x_k=kh$, $k=0,1,\dots,n$. 

One may propose one discretization of the solution of the form $u_k(t):= u(x_k,t)$, $k=0,1,\dots,n$. Then, the second derivative in \eqref{76} could be approximated using finite differences \cite{KIN}:

\begin{equation}\label{77}
\frac{\partial^2}{\partial x^2}\,u(x_k,t) = \frac{u(x_k-1,t)-2u(x_k,t) + u(x_{k+1},t)}{h^2}\,,
\end{equation}
while for the first spatial derivative, we have

\begin{equation}\label{78}
\frac{\partial}{\partial x}\,u(x_k,t) = \frac{u(x_k+1,t) - u(x_{k-1},t)}{2h}\,.
\end{equation}
with $k=1,2,\dots,n-1$, so that equation \eqref{76} takes the form

\begin{equation}\label{79}
\frac{d}{dt}\,U(t) =F(U(t))\,,
\end{equation} 
where $F$ is a square matrix of order $n-1$ and $U(t)= (u_1(t),u_2(t),\dots,u_{n-1}(t)$ with $u_k(t):= u(x_k,t)$, $k=1,2,\dots,n-1$ and initial values $U(0)=(u(x_1,0), u(x_2,0),\dots,u(x_{n-1},0))$ and initial value $u(x,0)=h(x)$.  If it were possible to write equation \eqref{79} in the form:

\begin{equation}\label{80}
\frac{d}{dt}\,U(t) = A(U(t)) \cdot U(t)\,,
\end{equation} 
then, we would be able to apply our method to find an approximate solution of  \eqref{80} on a given finite interval. This property is not fulfilled by the general convection diffusion equation \eqref{76}. However, it is satisfied by a particular choice of this type of parabolic equations: {\it the non-linear Burger's diffusion equation}, which is \cite{BUR,BA}:

\begin{equation}\label{81}
\frac{\partial}{\partial t}\, u(x,t) = \frac{\partial^2}{\partial x^2}\, u(x,t) - u(x,t)\, \frac{\partial}{\partial x}\, u(x,t)\,.
\end{equation}

We choose $[0,1]$ as the integration interval for the coordinate variable, so that $X=1$ as above. For the time variable, we use $t\in [0,1]$. This choice is made just for simplicity and as an example to implement our numerical calculations. We use the finite difference method, where the second derivative and first derivatives are replaced as in \eqref{77} and \eqref{78}, respectively.

Using \eqref{77} and \eqref{78} in \eqref{81}, we have for each of the nodes the following recurrence relation:

\begin{eqnarray}\label{82}
\frac{d}{dt}\, u(x_k,t) = \frac{1}{h^2} \left( \left( 1-\frac h2\, u(x_k,t) \right) u(x_{k-1},t) -2 u(x_k,t) + \left( 1 + \frac h2 \, u(x_k,t)  \right) u(x_{k+1},t) \right)\,,
\end{eqnarray}
for $k=1,2,\dots,n-1$. When written expressions \eqref{82} in matrix form, we obtain a matrix equation of the form \eqref{6} and, then, suitable for applying our method. 

In our numerical realization, we have used a small number of nodes, to begin with, say $n=5$. Then, we integrate on the time variable, on the interval $t\in [0,1]$, using $h_t:=1/m$, where $m$ is a given integer, as the distance between time nodes. Then, we compare our solution with the solution of \eqref{81} given by the sentence NDSolve provided by the Mathematica software, solution that we denote as $v(x,t)$. 

Then, we can estimate the error of our solution as compared with $v(x,t)$. This is given by

\begin{equation}\label{83}
{\rm error}:= \frac 1m \sum_{j=0}^m \sum_{k=1}^{n-1} (u(x_k,t_j)-v(x_k,t_j))^2\,,
\end{equation}
where $t_j:=jh_t$ for all $j=0,1,2,\dots,m$. To estimate the error, we have to give values to $m$. For $m=10$, the error is $1.07\, 10^{-4}$. A similar result can be obtained with $m=20$ or even higher, so that $m=10$ gives already a reasonable approximation. The solution for $t=1$ is nearly zero, which one may have expected taking into account that equation \eqref{80} describes a dissipative model. Thus, our approximate solution may be considered satisfactory also in this case. 

A few more words on the comparison of our solution and the solution using the Euler method, performed through NDSolve. First of all, we use the explicit Euler method, for which the local error is $o(h^2_t)$, through the option ``Method $\to$ ``ExplicitEuler'', ``StartingStepSize'' $\to$ $1/100$'', since for $h_t>0.01$ the result is unstable.  Once we have done the spatial discretization, we integrate \eqref{82} with respect to time. The instability often appears when one uses an explicit method of spatial and time discretization for parabolic PDE \cite{CLW}. This instability comes after the errors due to the arithmetic calculations and are amplified after time integration. Then, let us consider $n=5$, where $n$ is the number of spatial nodes, $0\le t \le 1$ and $m=100$, so that the time integration interval becomes $h_t=0.01$. Then, the error \eqref{83} is $3.40\, 10^{-4}$.

We may improve time integration using Euler mid-point integration. In this case, one uses the option ``Method $\to$ ``ExplicitMidpoint'', ``StartingStepSize'' $\to$ 1/10''. Choosing $m=10$, we obtain an error of $1.76\, 10^{-4}$, so that $h_t=0.1$. This error is of the order given by our method. Just recalling that the local error given by the mid-point Euler method is of the order of $o(h_t^3)$.

In addition, we may compare the precision of our method and both Euler methods mentioned above by the errors obtained using the third and forth order Adams method (see Chapter III in \cite{HNW}), which for $h_t=0.1$ are respectively given by $1.76\, 10^{-4}$ and $1.74\,10^{-4}$. This error has always been obtained using \eqref{83}, where $u(x,t)$ is our solution and $v(x,t)$ is the solution given by either Euler, Euler mid-point or Adams methods.  Nevertheless, the Adams method is a multi-step method while ours is one step method, which means that ours is more easily programmable.

}

}

\end{itemize}

\section{Concluding remarks}

In the present paper, we have generalized a one step integration method, which has been developed in the seventies of last century by Demidovich and some other authors \cite{DEM,KOT,HSU}. While Demidovich and others restrict themselves to the search for approximate solutions of linear systems of first order differential equations, albeit with variable coefficients, we propose a way to extend the ideas of the mentioned authors to non-linear systems.  The solutions we have obtained have a similar degree of precision than those proposed in \cite{DEM,KOT,HSU}. In our method, we obtain a sequence of approximate solutions on a finite interval of ordinary differential equations, and we have proven that this sequence converges uniformly to the exact solution. In order to obtain each approximate solution, we have divided the integration interval into subintervals of length $h$. The sequence of approximate solutions can be obtained after successive refinements of $h$. For each approximate solution, characterized by a value of $h$, we  have determined the local error up to $o(h^3)$. 

It is certainly true, that our one step matrix method does not improve the precision obtained with the fourth order Runge-Kutta method or other equivalent. However, the great advantage of the proposed method with any other is that its algorithm of construction, through the exponential matrix as described in Section 2, is much simpler than their competitors. Simplicity that is inherited from its antecesor the Demidovich method \cite{DEM,KOT,HSU}. This paper is somehow the continuation of previous research of the authors in the same field \cite{GL1,GLP}. 

We have added some examples of the application of the method, on where we have performed numerous numerical test on the precision of the method, which lies between second and third Taylor method. We have used the sofware Mathematica to implement these numerical tests.

\section*{Acknowledgements}

We acknowledge partial financial support to the Spanish MINECO, grant MTM2014-57129-C2-1-P, and the Junta de Castilla y Le\'on, grants  VA137G18, BU229P18 and VA057U16. We are grateful to Prof L.M. Nieto (Valladolid) for some useful suggestions.

\section{Appendix A: A conjecture relative to the Li\'enard equation.}

As is well known, the Li\'enard equation has the following form:

\begin{equation}\label{84}
y''(t)+f(y)\,y'(t) +g(y)=0\,.
\end{equation}

Let us assume that the function $g(y)$ is a product of some function, that we also call $g(y)$ for simplicity, and $y$, so that equation \eqref{67} takes the form:

\begin{equation}\label{85}
y''(t)+f(y)\,y'(t) +g(y)\,y(t)=0\,.
\end{equation}

This second order equation may be easily transformed into a two dimensional first order equation ($y(t)=y(t)$, $z(t)=y'(t)$):

\begin{equation}\label{86}
\left( \begin{array}{c} y' \\[2ex] z' \end{array} \right) = \left( \begin{array}{cc} 0 & 1 \\[2ex] -g(y) & -f(y)  \end{array} \right)  \left( \begin{array}{c} y \\[2ex] z \end{array} \right) = A \left( \begin{array}{c} y \\[2ex] z \end{array} \right)\,,
\end{equation}
where the meaning of the matrix $A$ is obvious. Its eigenvalues are

\begin{equation}\label{87}
\lambda_\pm(y) = -\frac12 \left(f(y) \pm \sqrt{f^2(y)-4g(y)}\,\right)\,.
\end{equation} 

The conjecture is the following: A sufficient condition for the solutions of \eqref{67} to be bounded for $t>0$ is that the following three properties hold simultaneously: i.) The discriminant in \eqref{60} be negative, i.e., $f^2(y)-4g(y)<0$; ii.) The function $f(y)$ be non-negative, $f(y)\ge 0$ and iii.) The function $g(y)$ is positive and smaller than one, $0< g(y)<1$. 

This conjecture is based in the form of equations \eqref{11} and \eqref{12} and we have tested it in several numerical experiments.

Two more comments in relation to the Li\'enard equation:

1.- The vector field associated to the equation is given by $J\equiv(z,-g(y)y-f(y)z)$, for which the divergence is ${\rm div}\, (J)=-f(y) \le 0$ if $f(y)\ge 0$. Therefore, if $f(y)$ is non-negative the divergence is always negative, so that the origin $(0,0)$ is an attractor. We conjecture that this attractor is also global.

2.- If we take $f(y)\equiv 0$, then \eqref{72} represents a Hamiltonian flow with Hamiltonian given by

\begin{equation}\label{88}
H(y,z)=\frac 12 \,z^2+V(y)\,,
\end{equation}
with

\begin{equation}\label{89}
V(y)=\int_0^y ug(u)\,du
\end{equation}

Under the assumption $g(y)>0$, the derivative $V'(y)=yg(y)$ is positive for $y>0$ and negative for $y<0$. Thus the only critical point is $y^*=0$, at this point $V'(0)=0$ and $V''(0)=g(0)>0$, so that all the orbits are closed, hence periodic.

\section{Appendix B: Source code to use our method in the van der Pole equation}

\begin{figure}
\centering
\includegraphics[width=1.4\textwidth]{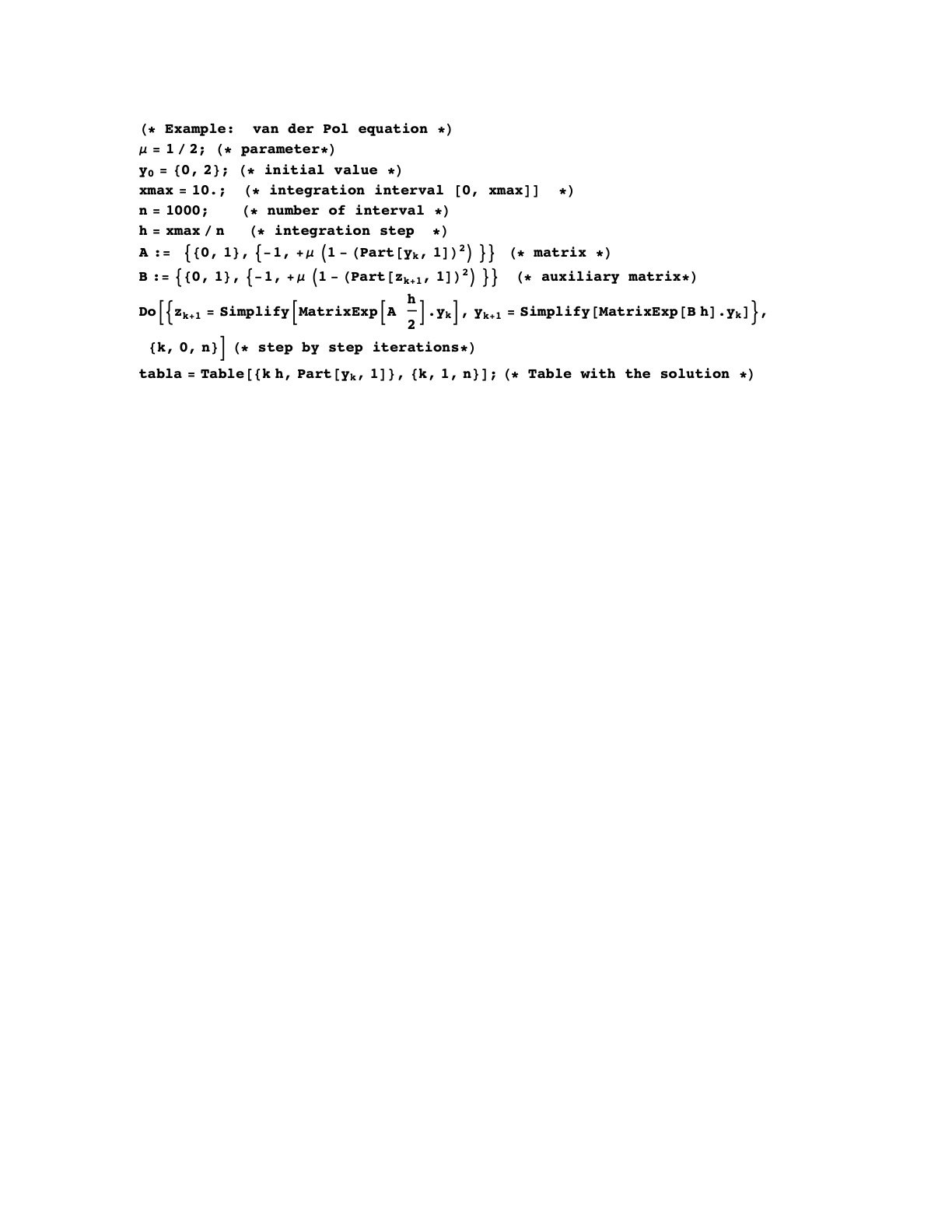}
\caption{\small Source Code used when applying the method to our example using the van der Pole equation.
\label{Figure5}}
\end{figure}

\vfill\eject


\begin{thebibliography}{99}

\bibitem{DEM}  Demidovich B.  {\it Lectures on the mathematical Theory of Stability}, Nauka, Moscow, 1967.

\bibitem{KOT}  Kotsis D.  The approximation of the characteristic multipliers of periodic differencial equations. {\it Alkalmaz. Mat. Lapok.}, 1977; 2: 269276.

\bibitem{HSU}  Hsu C.S. On approximating a general linear periodic system  {\it J. Math. Anal. Appl.} 1974; 45: 234251.

\bibitem{FAR}  Farkas M.  {\it Periodic Motions}. Springer, New York, 1994. 

\bibitem{VDP}  van der Pol B.  On relaxation-oscillations, {\it The London, Edinburgh and Dublin Phil. Mag. \& J. of Sci.} 1927;  2 (7): 978-992.

\bibitem{NEP} Nepomuceno E.P., Peixoto M.L.C., Martins S.A.M., Rodrigues Junior H.M., Perc M., Inconsitencies in numerical simulations of dynamical systems using interval arithmetic, {\it Int. J. Bifur. Chaos}  2018; 28 (04):  1850055. 

\bibitem{NEP1} Nepomuceno E.P., Guedes P.F.S., Barbosa A.M., Perc M., Repnik R.,  Soft computing simulations of chaotic systems, {\it Int. J. Bifur. Chaos} 2019; 29 (08): 1950112. 

\bibitem{DUF}  Duffing G. {\it Forced oscillations with variable natural frequency and their technical relevance} (German), Heft 41/42, Braunschweig: Vieweg; 1918, 1-134. 

\bibitem{LOR}  Lorenz E.N.,  Deterministic non-periodic flow, {\it J. Atmos. Sci.} 1963;  20 (2): 130-141. 

\bibitem{HAL}  Hale J.,  Kocak H., {\it Dynamics and Bifurcations}. Springer-Verlag; New York, 1991. 

\bibitem{LOT}  Lotka A.J. {\it Elements of Physical Biology}.  Williams and Wilkins; Philadelphia, USA, 1925.

\bibitem{VOL}  Volterra V. Variations and fluctuations of the number of individuals in animal species living together, in {\it Animal Ecology}, R.N. Chapman, Ed., McGraw-Hill, New York, 1931. 


\bibitem{BEL}  Bel\'endez A.,  Gimeno E.,  Fern\'andez E.,  M\'endez D.I.,  Alvarez M.L.  Accurate approximate solutions to non-linear oscillators in which the restoring force is inversely proportional to the dependent variable. {\it Phys. Scr.} 2008;  77: 065004 (2008).

\bibitem{EGG}  Egger H.,  Schmidt K.,  Shashkov V. Multistep and RungeKutta convolution quadrature methods for coupled dynamical systems, to be published in the {\it J. Comp. Appl. Math.}

\bibitem{KAL}  Kalogiratoua Z.,  Monovasilis Th.,  Psihoyios G.,  Simos T.E., Runge-Kutta type methods with special properties for the numerical integration of ordinary differential equations. {\it Phys. Rep.} 2014;  536: 75-146. 


\bibitem{MIC}   Mickens R. {\it Oscillations in Planar Dynamic Systems}. World Scientific,  London, 1996.



\bibitem{HAK}  Haken H.  Analogy between higher instabilities in fluids and lasers. {\it  Phys. Lett. A} 1975;  53: 77-78.

\bibitem{GOR}  Gorman M.,  Widmann P.J.,  Robbins K.A., Nonlinear dynamics of a convection loop: A quantitative comparison of experiment with theory. {\it  Physica D.} 1986;  19 (2): 255-267. 

\bibitem{COU}   Coumo K.M.,  Oppenheim A.V., Circuit implementation of synchronized chaos with applications to communications. {\it Phys. Rev. Lett.} 1993;  71 (1): 65-68. 

\bibitem{MIS}  Mishra A.,  Sanghi S.  A study of the asymmetric Malkus waterwheel: The biased Lorenz equations. {\it Chaos} 2006;  16 (1): 013114. 






\bibitem{BAR}  Bario R.  Performance of the Taylor series method for ODEs/DAEs. {\it Appl. Math. Comp.} 2005; 163: 525-545.

\bibitem{LG}  Lara L.P.,  Gadella M.  An approximation to solutions of linear ODE by cubic interpolation{\it  Comp. Math. Appl.} 2008;  56: 1488-1495. 

\bibitem{KK}  Kermack W.O.,  McKendrick A.G.  A contribution to the Mathematic Theory of Epidemics. {\it Proceedings of the Royal Society A} 1927;  115 (772): 700-721. 

\bibitem{STR} S. Strogatz. {\it Non-linear Dynamics and Chaos}.  Addison-Wesley, New York, 1994.

\bibitem{TC}  Torralba-Rodr\'iguez O.,  Conde-Guti\'errez R.A.,  Hern\'andez-Javier A.L.  Modeling and prediction of COVID-19 in M\'exico applying mathematical and computational models. {\it Chaos, Solitons and Fractals} 2020; 138: 1009946 (2020). 



\bibitem{CI}  Chicone C. {\it Ordinary Differential Equations with Applications}. Springer, New York, 1999. 

\bibitem{MUR}  Murray J.D. {\it Mathematical Biology I: An Introduction}. Springer, Berlin and New York, 2003.

\bibitem{LV}  Llibre J.,  Valls C., Global analytic first integrals for the real planar Lotka-Volterra system. {\it J. Math. Phys.} 2007;  48 (3): 033507. 

\bibitem{HWX} Huabing Jia, Wei Xu, Xiaoshan Zhao, Zhanguo Li, Separation of variables and exact solutions to nonlinear diffusion equations with x-dependent convection and absorption. {\it J. Math. Annal. Appl.} 2008; 339: 982-995. 

\bibitem{KIN} Kincaid D, Cheney W, {\it Numerical analysis, Mathematics of Scientific Computing}, AMS, Providence, Rhode Island, USA, Third Edition 2002. 

\bibitem{BUR} Burger J.M., {\it A mathematical Model Illustrating the Theory of Turbulence}, Academic Press, New York, 1948.

\bibitem{BA} Biazar J., Aminikhah H., Exact numerical solutions for non-linear Burger's equation by VIM. {\it Math. Comp. Mod.} 2009; 49: 1394-1400. 

\bibitem{CLW} Carnahan B., Lutter H.A., Wilkes J.O., {\it Applied Numerical Methods}. Wiley, New York, 1969.

\bibitem{HNW} Hairer E., N{\o}rsett S.P., Wanner G., {\it Solving Ordinary Differential equations I: Nonstiff Problems}. Springer, Berlin, Heidelberg, 1993. 

\bibitem{GL1}  Gadella M.,  Lara L.P.  On the determination of approximate periodic solutions of some non-linear ODE. {\it Appl. Math. Comp.} 2012;   18 (10): 6038-6044.

\bibitem{GLP}  Gadella M.,  Lara L.P.,  Pronko G.P.  Iterative solution of some nonlinear differential equations. {\it Appl. Math. Comp.} 2011;  217 (22): 9480-9487. 






\end{thebibliography}
\end{document}